\documentclass[11pt]{article}

\usepackage[utf8]{inputenc}
\usepackage{amsmath,amsthm}
\usepackage{amssymb}
\usepackage{geometry}
\usepackage{bbm}

\def\er{{\mathbb R}}
\def\zet{{\mathbb Z}}
\def\Ex{{\mathbb E}}
\def\Pr{{\mathbb P}}
\def\de{\mathrm{d}}
\def\Var{\mathrm{Var}}

\def\conv{\mathrm{conv}}
\def\dim{\mathrm{dim}}
\def\vol{\mathrm{vol}}

\newtheorem{thm}{Theorem}
\newtheorem{prop}[thm]{Proposition}
\newtheorem{lem}[thm]{Lemma}

\theoremstyle{remark}
\newtheorem{rem}[thm]{Remark}

\title{Bounding suprema of canonical processes via convex hull
\thanks{Supported by the National Science Centre, Poland grant 2015/18/A/ST1/00553}}
\author{Rafał Latała}
\date{}

\begin{document}

\maketitle

\begin{abstract}
We discuss the method of bounding suprema of canonical processes based on the inclusion of their index set into a convex hull of a well-controlled set of points. While the upper bound is immediate, the reverse estimate was established to date 
only for a narrow class of regular stochastic processes. We show that for specific index sets, including arbitrary ellipsoids, regularity assumptions may be substantially weakened.
\end{abstract}

\section{Formulation of the problem}

Let $X=(X_1,\ldots,X_n)$ be a centered random vector with independent coordinates. To simplify the notation we
will write 
\[
X_t=\langle t,X\rangle=\sum_{i}t_iX_i\quad \mbox{ for }t=(t_1,\ldots,t_n)\in\er^n.
\]
Our aim is to estimate the expected value of the supremum of the process $(X_t)_{t\in T}$, i.e. the quantity
\[
b_X(T):=\Ex\sup_{t\in T}X_t, \quad T\subset \er^n \mbox{ nonempty bounded}.
\]
There is a long line of research devoted to bounding $b_X(T)$ via the chaining method (cf. the monograph \cite{Tabook}). However chaining methods do not work well for heavy-tailed random variables. In this
paper we will investigate another approach based on the convex hull method.

First let us discuss an easy upper bound. Suppose that there exists $t_0,t_1,\ldots\in\er^n$ such that
\begin{equation}
\label{eq:subsetch}
T-t_0\subset\overline{\mathrm{conv}}\{\pm t_i\colon\ i\geq 1\}  
\end{equation}
then for any $u>0$,
\[
\Ex\sup_{t\in T}X_t=\Ex\sup_{t\in T}X_{t-t_0}
\leq \Ex\sup_{i\geq 1}|X_{t_i}|
\leq u+\sum_{i\geq 1}\Ex|X_{t_i}|I_{\{|X_{t_i}|\geq u\}}.
\]
Indeed the equality above follows since $X_{t-t_0}=X_t-X_{t_0}$ and $\Ex X_{t_0}=0$ and all inequalities are pretty obvious.
To make the notation more compact let us define for nonempty countable sets $S\subset \er^n$
\[
M_X(S)=\inf_{u>0}\Big[u+\sum_{t\in S}\Ex|X_{t}|I_{\{|X_{t}|\geq u\}}\Big],\quad
\widetilde{M}_X(S)=\inf\Big\{m>0\colon\ \sum_{t\in S}\Ex|X_{t}|I_{\{|X_{t}|\geq m\}}\leq m\Big\}.
\]

It is easy to observe that
\begin{equation}
\label{eq:compmM}
\widetilde{M}_X(S)\leq M_X(S)\leq 2\widetilde{M}_X(S).
\end{equation}
To see the lower bound let us fix $u>0$ and set $m=u+\sum_{t\in S}\Ex|X_{t}|I_{\{|X_{t}|\geq u\}}$
then 
\[
\sum_{t\in S}\Ex|X_{t}|I_{\{|X_{t}|\geq m\}}\leq \sum_{t\in S}\Ex|X_{t}|I_{\{|X_{t}|\geq u\}}\leq m,
\]
so $\widetilde{M}_X(s)\leq m$. For the upper bound it is enough to observe that for $u>\widetilde{M}_X(S)$ we have 
$\sum_{t\in S}\Ex|X_{t}|I_{\{|X_{t}|\geq u\}}\leq u$.

We have thus shown that
\begin{equation}
\label{eq:trivialupper}
b_X(T)\leq M_X(S)\leq 2\widetilde{M}_X(S)\quad \mbox{if}\quad T-t_0\subset \overline{\conv}(S\cup -S). 
\end{equation}

\begin{rem}
The presented proof of \eqref{eq:trivialupper} did not use independence of coordinates of $X$, the only required property
is mean zero.
\end{rem}

\noindent
\textbf{Main question.} \textit{When can we reverse bound \eqref{eq:trivialupper} -- what should be assumed about variables $X_i$
(and the set $T$) in order that
\begin{equation}
\label{eq:revest}
T-t_0\subset \overline{\mathrm{conv}}(S\cup -S)\quad \mbox{and}\quad  M_X(S)\lesssim \Ex\sup_{t\in T}X_t
\end{equation}
for some $t_0\in \er^n$ and nonempty countable set $S\subset \er^n$?}

\begin{rem}
It is not hard to show (see Section \ref{sec:B1} below) that $M_X(S)\sim \Ex \max_{i}|X_{t_i}|=b_X(S\cup -S)$
if $S=\{t_1,\ldots,t_k\}$ and variables $(X_{t_i})_{i}$ are independent. Thus our main question asks
whether the parameter $b_X(T)$ may be explained by enclosing a translation of $T$ into the convex hull
of points $\pm t_i$ for which variables $X_{t_i}$ behave as though they are independent.  
\end{rem}

\begin{rem}
The main question is related to Talagrand conjectures about suprema of positive selector processes, c.f. \cite[Section 13.1]{Tabook}, i.e. the case when $T\subset \er_+^n$ and $\Pr(X_i\in \{0,1\})=1$. Talagrand investigates possibility of enclosing $T$ into a solid convex hull, which is  bigger than the convex hull. On the other hand we think that in our question some regularity conditions on variables $X_i$ is needed (such as $4+\delta$ moment condition \eqref{eq:4+mom}, which is clearly not satisfied for nontrivial classes of selector processes). 
\end{rem}

\begin{rem}
i) In the one dimensional case if $a=\inf T$, $b=\sup T$, then 
$T\subset [a,b]=\frac{a+b}{2}+\conv\{\frac{a-b}{2},\frac{b-a}{2}\}$.
Hence
\[
b_{X_1}(T)=\Ex\max\{aX_1,bX_1\}=\frac{a+b}{2}\Ex X_1+\Ex\Big|\frac{b-a}{2}X_1\Big|=\frac{b-a}{2}\Ex |X_1|
\geq \widetilde{M}_{X_1}\Big(\Big\{\frac{b-a}{2}\Big\}\Big),
\]
so this case is trivial. Thus in the sequel it is enough to consider $n\geq 2$.
\smallskip
\\
ii) The set $V:=\overline{\conv}(S\cup -S)$ is convex and origin-symmetric. Hence
if $T=-T$ and $T-t_0\subset V$ then $T+t_0=-(-T-t_0)=-(T-t_0)\subset V$
and $T\subset\conv((T-t_0)\cup (T+t_0))\subset V$. Thus for symmetric sets it is enough to consider only $t_0=0$.
\smallskip
\\
iii) Observe that $b_X(\conv(T))=b_X(T)$ and $T-t_0$ is a subset of a convex set if and only if $\conv(T)-t_0$ is a subset of this set.
Moreover, if $T-T\subset V$ then $T-t_0\subset V$ for any $t_0\in V$ and
$b_X(T-T)=b_X(T)+b_X(-T)=b_X(t)+b_{-X}(T)$. So if $X$ is symmetric it is enough to consider symmetric convex sets $T$.
\end{rem}

\noindent
\textbf{Notation.} Letters $c$, $C$ will denote absolute constants which value may differ at each occurence. For two nonnegative functions $f$ and $g$ we write $f\gtrsim g$ (or $g\lesssim f$) if  $g\leq Cf$. Notation $f\sim g$ means that 
$f\gtrsim g$ and $g\gtrsim f$. We write $c(\alpha)$, $C(\alpha)$ for constants depending only on a parameter $\alpha$ and 
define accordingly relations $\gtrsim_\alpha$, $\lesssim_\alpha$,  $\sim_\alpha$. 

\medskip

\noindent
\textbf{Organization of the paper.} In Section \ref{sec:regular} we present another quantity $m_X(S)$, 
defined via $L_p$-norms of $(X_t)_{t\in S}$, and show that for regular variables $X_i$ it is equivalent to $M_X(S)$. We also discuss there the relation of the convex hull method to the chaining functionals. In Section \ref{sec:B1} we show that for $T=B_1^n$ the bound \eqref{eq:trivialupper} may be reversed for arbitrary independent $X_1,\ldots,X_n$ and $S=\{e_1,\ldots,e_n\}$.
Section \ref{sec:B2} is devoted to the study of ellipsoids. First we show that for $T=B_2^n$ and symmetric $p$-stable random
variables, $1<p<2$, one cannot reverse \eqref{eq:trivialupper}. Then we prove that under $4+\delta$ moment condition our main
question have the affirmative answer for $T=B_2^n$ and more general case of ellipsoids. We extend this result to the case of
linear images of $B_q^n$-balls, $q\geq 2$ in Section \ref{sec:Bq}. We conclude by discussing some open questions in the last section.

\section{Regular growth of moments.}
\label{sec:regular} 

In this section we consider variables with regularly growing moments in a sense that
\begin{equation}
\label{eq:regmom}
\|X_i\|_{2p}\leq \alpha\|X_i\|_p<\infty\quad \mbox{ for }p\geq 1,
\end{equation}
where $\|X\|_p=(\Ex |X|^p)^{1/p}$.

For such variables we will prove that there is alternate quantity equivalent to $M_X(S)$, namely 
\[
m_X(S):=\inf\sup_{i}\|X_{t_i}\|_{\log(e+i)}.
\]
where the infimum runs over all numerations of $S=\{t_i\colon\ 1\leq i\leq N\}$, $N\leq \infty$. 

It is not hard to check (cf. Lemma 4.1 in \cite{LaSt}) that \eqref{eq:regmom} yields
\begin{equation}
\label{eq:regmom2}
\|X_t\|_{2p}\leq C_0(\alpha)\|X_t\|_p\quad \mbox{ for }p\geq 1 
\end{equation}
and as a consequence we have for $p>0$,
\begin{equation}
\label{eq:tailmomreg}
\Pr(|X_t|\geq e\|X_t\|_p)\leq e^{-p},\quad \Pr(|X_t|\geq c_1(\alpha)\|X_t\|_p)\geq \min\{c_2(\alpha),e^{-p}\},
\end{equation}
where the first bound follows by Chebyshev's inequality and the second one by the Paley-Zygmund inequality.

\begin{prop}
Suppose that $X_i$ are independent r.v's satisfying condition \eqref{eq:regmom}. Then
$M_X(S)\sim_\alpha m_X(S)$.
\end{prop}

\begin{proof}
Let $S=\{t_i\colon\ 1\leq i\leq N\}$ and $m:=\sup_{i}\|X_{t_i}\|_{\log(e+i)}$. Then for $u>1$,
\[
\sum_{s\in S}\Pr(|X_s|\geq um)
\leq\sum_{i=1}^N\Pr(|X_{t_i}|\geq u\|X_{t_i}\|_{\log(e+i)})\leq
\sum_{i=1}^N u^{-\log(e+i)}.
\]
Therefore
\begin{align*}
\sum_{s\in S}\Ex|X_s|I_{\{|X_s|\geq e^2m\}}
&=\sum_{s\in S}\Big(e^2m\Pr(|X_s|\geq e^2m)+m\int_{e^2}^\infty\Pr(|X_s|\geq um)\de u\Big)
\\
&\leq m\sum_{i=1}^N\Big(e^{2-2\log(e+i)}+\int_{e^2}^\infty u^{-\log(e+i)}\de u\Big)
\\
&\leq m\sum_{i=1}^N\Big((e+i)^{-2}\Big(e^2+\frac{1}{\log(e+i)-1}\Big)\Big)\leq 100m,
\end{align*}
which shows that $M_X(S)\leq 100m_X(S)$ (this bound does not use neither regularity neither independence of $X_i$).

To establish the reverse inequality let us take any $m>2M_X(S)\geq\widetilde{M}_X(S)$ and enumerate elements of
$S$ as $t_1,t_2,\ldots$ in such a way that that $i\rightarrow \Pr(|X_{t_i}|\geq m)$ is nonincreasing. By the definition
of $\widetilde{M}_X(S)$ we have 
\[
\sum_{i=1}^N\Pr(|X_{t_i}|\geq m)\leq \frac{1}{m}\sum_{i=1}^N\Ex|X_{t_i}|I_{\{|X_{t_i}|\geq m\}}\leq 1.
\]
In particular it means that $\Pr(|X_{t_i}|\geq m)\leq 1/i$. By \eqref{eq:tailmomreg} this yields that for 
$i>1/c_2(\alpha)$
$\|X_{t_i}\|_{\log(i)}\leq m/c_1(\alpha)$. Since $\log(e+i)/\log(i)\leq 2$ for $i\geq 3$ we have
$\|X_{t_i}\|_{\log(e+i)}\leq C(\alpha)m$ for large $i$. For $i\leq \max\{3,1/c_2(\alpha)\}$ it is enough to observe
that $\log(e+i)\leq 2^{k(\alpha)}$, so 
\[
\|X_{t_i}\|_{\log(e+i)}\leq C_0(\alpha)^{k(\alpha)}\Ex|X_{t_i}|\leq  C_0(\alpha)^{k(\alpha)}M_X(S).
\]
This shows that $\|X_{t_i}\|_{\log(e+i)}\lesssim_\alpha m$ for all $i$ and therefore $m_X(S)\lesssim_\alpha M_X(S)$.

\end{proof}

\subsection{$\gamma_X$-functional}
\label{sub:gammaX}

The famous Fernique-Talagrand  theorem \cite{Fe,Ta_reg} states that suprema of Gaussian processes may be estimated in geometrical
terms by $\gamma_2$-functional. This result was extended in several directions. One of them is based on  the so-called
$\gamma_X$ functional.

For a nonempty subset $T\subset \er^n$ we define
\[
\gamma_X(T):=\inf\sup_{t\in T}\sum_{n=0}^\infty \Delta_{n,X}(A_n(t)),
\]
where the infimum runs over all increasing sequences of partitions $({\cal A}_n)_{n\geq 0}$ of $T$ such that
${\cal A}_{0}=\{T\}$ and $|{\cal A}_n|\leq N_n:=2^{2^n}$ for $n\geq 1$, $A_n(t)$ is the unique element of ${\cal A}_n$ which contains $t$
and $\Delta_{n,X}(A)$ denotes the diameter of $A$ with respect to the distance $d_n(s,t):=\|X_s-X_t\|_{2^n}$.

It is not hard to check that $b_X(T)\lesssim \gamma_X(T)$. The reverse bound was discussed in
\cite{LaTk}, where it was shown that it holds (with constants depending on $\beta$ and $\lambda$)
if 
\begin{equation}
\label{eq:regmomgr}
\|X_i\|_p\leq \beta\frac{p}{q}\|X_i\|_q \mbox{ and }\|X_i\|_{\lambda p}\geq 2\|X_i\|_p
\mbox{ for all $i$ and $p\geq q\geq 2$.}
\end{equation}
Moreover the condition $\|X_i\|_p\leq \beta\frac{p}{q}\|X_i\|_q$ is necessary in the i.i.d. case if the estimate
$\gamma_X(T)\leq C b_X(T)$ holds with a constant independent on $n$ and $T\subset \er^n$. 
  
The next result may be easily deduced from the proof of \cite[Corollary 2.7]{LaTk}, but we provide its proof for the sake
of completeness. 

\begin{prop}
\label{prop:gammaX2MX}
Let $X_i$ be independent and satisfy condition \eqref{eq:regmom} and let $T$ be a nonempty subset of $\er^n$ such that $\gamma_X(T)<\infty$. Then there exists set $S\subset \er^n$ 
such that for any $t_0\in T$, $T-t_0\subset T-T\subset \overline{\conv}(S\cup -S)$ and $M_X(S)\lesssim m_X(S)\lesssim_\alpha \gamma_X(T)$. 
\end{prop}
 
\begin{proof}
Wlog (since it is only a matter of rescaling) we may assume that $\Ex X_i^2=1$.
 
By the definition of $\gamma_X(T)$ we may find an increasing sequence of partitions $({\cal A}_n)$ such that
${\cal A}_0=\{T\}$, $|{\cal A}_j|\leq N_j$ for $j\geq 1$ and 
\begin{equation}
\label{eq:part}
\sup_{t\in T}\sum_{n=0}^\infty \Delta_{n,X}(A_n(t))\leq 2\gamma_X(T).
\end{equation}
For any $A\in {\cal A}_n$ let us choose a point $\pi_n(A)\in A$ and set $\pi_n(t):=\pi_n(A_n(t))$. 

Let $M_n:=\sum_{j=0}^{n} N_j$ for $n=0,1,\ldots$ (we put $N_0:=1$). 
Then $\log(M_n+2)\leq 2^{n+1}$.
Notice that there are $|{\cal A}_n|\leq N_n$ points of the form $\pi_n(t)-\pi_{n-1}(t)$, $t\in T$.
So we may define
$s_{k}$, $M_{n-1}\leq k< M_n$, $n=1,2,\ldots$ as some rearrangement (with repetition if $|{\cal A}_n|<N_n$) 
of points of the form $(\pi_n(t)-\pi_{n-1}(t))/\|X_{\pi_n(t)}-X_{\pi_{n-1}(t)}\|_{2^{n+1}}$, $t\in T$.
Then $\|X_{s_k}\|_{\log(k+e)}\leq 1$ for all $k\geq 1$.

Observe that 
\[
\|t-\pi_n(t)\|_2=\|X_t-X_{\pi_n(t)}\|_{2}\leq \Delta_{2,X}(A_n(t))\leq \Delta_{n,X}(A_n(t))\rightarrow 0
\quad \mbox{ for }n\rightarrow\infty.
\]
For any $s,t\in T$ we have $\pi_0(s)=\pi_0(t)$ and thus
\[
s-t=\lim_{n\rightarrow\infty}(\pi_n(s)-\pi_n(t))=
\lim_{n\rightarrow\infty}\left( \sum_{k=1}^n(\pi_k(s)-\pi_{k-1}(s))-\sum_{k=1}^n(\pi_k(t)-\pi_{k-1}(t))\right).
\]
This shows that
\[
T-T\subset R\ \overline{\mathrm{conv}}\{\pm s_k\colon\ k\geq 1\}, 
\]
where
\begin{align*}
R&:=2\sup_{t\in T}\sum_{n=1}^{\infty}d_{n+1}(\pi_n(t),\pi_{n-1}(t))
\leq 2\sup_{t\in T}\sum_{n=1}^{\infty}\Delta_{n+1,X}(A_{n-1}(t))
\\
&\leq  C(\alpha)\sup_{t\in T}\sum_{n=1}^{\infty}\Delta_{n-1,X}(A_{n-1}(t))
\leq 2C(\alpha)\gamma_X(T),
\end{align*}
where the second inequality follows by \eqref{eq:regmom2}.
Thus it is enough to define $S:=\{Rs_k\colon\ k\geq 1\}$.
\end{proof}

\begin{rem}
Proposition \ref{prop:gammaX2MX} together with the equivalence $b_X(T)\sim_{\alpha,\lambda}\gamma_X(T)$ shows that  the main question
has the affirmative answer for any bounded nonempty set $T$ if symmetric random variables $X_i$ satisfy moment bounds \eqref{eq:regmom}.
We strongly believe that the condition $\|X_i\|_{\lambda p}\geq 2\|X_i\|_p$ is not necessary -- equivalence of $b_X(T)$ and the convex hull bound was established in the case of symmetric Bernoulli r.v's ($\Pr(X_i=\pm 1)=1/2$) in \cite[Corollary 1.2]{BeLa}. However to treat the general case of r.v's satisfying only the condition $\|X_i\|_{p}\leq \beta\frac{p}{q}\|X_i\|_q$ 
one should most likely combine $\gamma_X$ functional with a suitable decomposition of the process $(X_t)_{t\in T}$, as was done for Bernoulli processes.
\end{rem}

\section{Toy case: $\ell_1$-Ball} 
\label{sec:B1}

Let us now consider a simple case of $T=B_1^n=\{t\in\er^n\colon\ \|t\|_1\leq 1\}$.
Let 
\[
u_0:=\inf\Big\{u>0\colon\ \Pr\big(\max_i|X_i|\geq u\big)\leq \frac{1}{2}\Big\}.
\]
Since
\[
\Pr\big(\max_i|X_i|\geq u\big)\geq \frac{1}{2}\min\Big\{1,\sum_{i}\Pr(|X_i|\geq u)\Big\}
\]
we get
\begin{align*}
\Ex\sup_{t\in B_1^n}X_t
&=\Ex\max_{1\leq i\leq n}|X_i|=\int_{0}^\infty\Pr\big(\max_{1\leq i\leq n}|X_i|\geq u\big)\de u
\geq \frac{1}{2}u_0+\int_{u_0}^\infty\frac{1}{2}\sum_{i=1}^n\Pr(|X_i|\geq u)\de u
\\
&= \frac{1}{2}u_0+\frac{1}{2}\sum_{i=1}^n\int_{u_0}^\infty\Pr(|X_i|\geq u)\de u=
\frac{1}{2}u_0+\frac{1}{2}\sum_{i=1}^n\Ex(|X_i|-u_0)_{+}.
\end{align*}
Therefore
\[
2u_0+\sum_{i=1}^n\Ex|X_i|I_{\{|X_i|\geq 2u_0\}}\leq 2u_0+2\sum_{i=1}^n\Ex(|X_i|-u_0)_{+}\leq 4\Ex\sup_{t\in B_1^n}X_t,
\]
so that $M_X(\{e_i\colon\ i\leq n\})\leq 4\Ex\sup_{t\in B_1^n}X_t$, where $(e_i)_{i\leq n}$ is the canonical basis of $\er^n$.
Since $B_1^n\subset \conv\{\pm e_1,\ldots,\pm e_n\}$ we get the affirmative answer to the main question for $T=B_1^n$.

\begin{prop}
If $T=B_1^n$ then estimate \eqref{eq:revest} holds for arbitrary independent integrable r.v's $X_1,\ldots,X_n$ with $S=\{e_1,\ldots, e_n\}$ and $t_0=0$.
\end{prop}

\section{Case II. Euclidean balls} 
\label{sec:B2}

Now we move to the case $T=B_2^n$. Then $\sup_{t\in T}\langle t,x\rangle=|x|$, where $|x|=\|x\|_2$ is the Euclidean norm
of $x\in \er^n$.

\subsection{Counterexample}
In this subsection $X=(X_1,X_2,\ldots,X_n)$, where $X_k$ have symmetric $p$-stable distribution with characteristic function
$\varphi_{X_k}(t)=\exp(-|t|^p)$ and $p\in (1,2)$. We will assume for convenience that $n$ is even. Let $G$ be a canonical $n$-dimensional Gaussian vector, independent of $X$. Then
\begin{align*}
\Ex|X|
&=\Ex_X\Ex_G\sqrt{\frac{\pi}{2}}|\langle X,G\rangle|=\sqrt{\frac{\pi}{2}}\Ex_G\Ex_X|\langle X,G\rangle|
=\sqrt{\frac{\pi}{2}}\Ex_G\|G\|_p\Ex|X_1|
\\
&\sim_p \Ex\|G\|_p\sim (\Ex\|G\|_p^p)^{1/p}\sim n^{1/p}.
\end{align*}

Observe also that for $u>0$, $\Pr(|X_1|\geq u)\sim_p \min\{1,u^{-p}\}$, so
\[
\Ex|X_1|I_{\{|X_1|\geq u\}}
\sim_p u\min\{1,u^{-p}\}+\int_u^\infty \min\{1,v^{-p}\}\de v\sim_p \min\{1,u^{1-p}\},\quad u>0
\]
and
\[
\Ex|X_t|I_{\{|X_t|\geq u\}}=\|t\|_p\Ex|X_1|I_{\{|X_1|\geq u/\|t\|_p\}}\sim_p 
\min\{\|t\|_p,u^{1-p}\|t\|_p^{p}\}, \quad u>0,\ t\in \er^n.
\]
Hence 
\[
\sum_{t\in S}\|t\|_p^p\lesssim_p u^p \quad \mbox{ for }u>\widetilde{M}_X(S).
\]

Suppose that $B_2^n\subset \overline{\conv}(S\cup -S)$ and  $M_X(S)\sim \widetilde{M}_X(S)<\infty$. 
We may then enumerate elements of $S$ as 
$(t_k)_{k=1}^N$, $N\leq \infty$ in such a way that $(\|t_k\|_p)_{k=1}^N$ is nonincreasing. Obviously $N\geq n$ (otherwise $\conv(S\cup -S)$ would have empty interior). Take $u>\widetilde{M}_X(S)$ and set 
$E:=\mathrm{span}(\{t_k\colon\ k\leq n/2\})$. 
Then $\|t_k\|_p^p\leq C_pu^p/n$ for $k> n/2$. Thus
\[
B_2^n\subset \overline{\conv}(S\cup -S)\subset E+\overline{\conv}(\{\pm t_k\colon\ k> n/2\})
\subset E+\Big(\frac{C_p}{n}\Big)^{1/p}uB_p^n.
\]
Let $F=E^\perp$ and $P_{F}$ denotes the ortogonal projection of $\er^n$ onto the space $F$.
Then $\dim F=\dim E=n/2$ and
\[
B_2^n\cap F=P_F(B_2^n)\subset \Big(\frac{C_p}{n}\Big)^{1/p}uP_F(B_p^n).
\]
In particular
\[
n^{-1/2}\sim \vol_{n/2}^{2/n}(B_2^n\cap F)\leq \Big(\frac{C_p}{n}\Big)^{1/p}u
 \vol_{n/2}^{2/n}(P_F(B_p^n)).
\]
By the Rogers-Shephard inequality \cite{RS} and inclusion $B_2^n\subset n^{1/p-1/2}B_p^n$ we have
\[
\vol_{n/2}(P_F(B_p^n))\leq \binom{n}{n/2}\frac{\vol_n(B_p^n)}{\vol_{n/2}(B_p^n\cap E)}
\leq 2^n\frac{\vol_n(B_p^n)}{\vol_{n/2}(n^{1/2-1/p}B_2^n\cap E)}\leq (Cn^{-1/p})^{n/2}.
\]
This shows that $u\gtrsim_p n^{2/p-1/2}$. Thus $M_X(S)\gtrsim_p n^{2/p-1/2}\gg n^{1/p}\sim_p b_X(B_2^n)$ and our question has a negative answer in this case.

\subsection{$4+\delta$ moment condition}

In this part we establish positive answer to the main question in the case $T=B_2^n$ under the following $4+\delta$ moment condition
\begin{equation}
\label{eq:4+mom}
\exists_{r\in (4,8],\lambda<\infty}\ (\Ex X_i^r)^{1/r}\leq \lambda (\Ex X_i^2)^{1/2}<\infty \quad i=1,\ldots,n.
\end{equation}

The restriction $r\leq 8$ is just for convenience. The following easy consequence of \eqref{eq:4+mom} will be helpful in the sequel.

\begin{lem}
\label{lem:estmom4+}
Suppose that $X_1,\ldots,X_n$ are independent mean zero r.v's satisfying condition \eqref{eq:4+mom}. Then for any $1\leq p\leq r$, 
\begin{equation}
\label{eq:complinear}
\Big\|\sum_{i=1}^nu_iX_i\Big\|_p\sim_\lambda \Big\|\sum_{i=1}^nu_iX_i\Big\|_2=\Big(\sum_{i=1}^n u_i^2\Ex X_i^2\Big)^{1/2}
\end{equation}
and
\begin{equation}
\label{eq:compquadr}
\Big\|\sum_{1\leq i<j\leq n}u_{ij}X_iX_j\Big\|_p\sim_\lambda \Big\|\sum_{1\leq i<j\leq n}u_{ij}X_iX_j\Big\|_2
=\Big(\sum_{1\leq i<j\leq n}u_{ij}^2\Ex X_i^2 \Ex X_j^2\Big)^{1/2}.
\end{equation}
\end{lem}

\begin{proof}
Since it is only a matter of scaling wlog we may and will assume that $\Ex X_i^2=1$ for all $i$.

Rosenthal's inequality \cite{Ro} gives for $2\leq p\leq r$ (recall that $r\in (4,8]$, so constants below do not depend on $r$)
\begin{align*}
\Big\|\sum_{i=1}^nu_iX_i\Big\|_p
&\sim \Big(\sum_i \Ex|u_iX_i|^2\Big)^{1/2}+\Big(\sum_i \Ex|u_iX_i|^p\Big)^{1/p}
\sim_\lambda \Big(\sum_i u_i^2\Big)^{1/2}+\Big(\sum_i |u_i|^p\Big)^{1/p}
\\
&\sim \Big(\sum_i u_i^2\Big)^{1/2}.
\end{align*}
 To estimate $\|S\|_p$ for  $1\leq p\leq 2$ and $S= \sum_{i=1}^nu_iX_i$ it is enough to note that
$\|S\|_1\leq \|S\|_p\leq \|S\|_2$ and $\|S\|_2\leq \|S\|_4^{1/3}\|S\|_1^{2/3}\sim_\lambda \|S\|_2^{1/3}\|S\|_1^{2/3}$, so $\|S\|_p\sim \|S\|_2$.

To prove the last part of the assertion we will use the hypercontractive method. Observe that for a real number $u$
there exists $\theta\in [0,1]$ such that
\begin{align*}
(1+u)^r&\leq \Big(1+ru+\frac{r(r-1)}{2}(1+\theta u)^{r-2}u^2\Big)I_{\{|u|<1\}}+(2|u|)^rI_{\{|u|\geq 1\}}
\\
&\leq 1+ru+r^22^{r-3}u^2+2^r|u|^r.
\end{align*}
Hence (note that $\lambda\geq 1$, $\Ex X_i=0$, $\Ex X_i^2=1$ and $\Ex|X_i|^r\leq \lambda^r$)
\[
\Ex\Big(1+\frac{1}{32\lambda}uX_i\Big)^r\leq 1+r^22^{r-3}\frac{u^2}{1024}+2^{-4r}|u|^r\leq 1+\frac{ru^2}{4}+\frac{|u|^r}{2}
\leq 1+\max\Big\{\frac{r}{2}u^2,|u|^r\Big\}.
\]
Since
\[
(\Ex(1+uX_i)^2)^{r/2}=(1+u^2)^{r/2}\geq 1+\max\Big\{\frac{r}{2}u^2,|u|^r\Big\}
\]
we get $\|1+\frac{1}{32\lambda}uX_i\|_r\leq \|1+uX_i\|_2$ for any $u\in \er$ and the hypercontractivity method
(cf. \cite[Theorem 6.5.2]{KW}) yields \eqref{eq:compquadr} for $p=r$. The case $1\leq p\leq r$ may be obtained 
in the same way as in the proof of \eqref{eq:complinear}.
\end{proof}


Observe that \eqref{eq:4+mom} implies that $\Var(X_i^2)\leq (\lambda^4-1)(\Ex X_i^2)^2$, so
$\Var(|X|^2)\leq \sum_{i}(\lambda^4-1)(\Ex X_i^2)^2\leq (\lambda^4-1)(\Ex|X|^2)^2$. This yields that
$\Ex|X|^4\leq \lambda^4(\Ex|X|^2)^2$ and $(\Ex|X|^2)^{1/2}\leq \lambda^2\Ex|X|$.

The next fact is pretty standard, we prove it for completeness.

\begin{lem}
\label{lem:net}
For any $k$ there exists $T\subset B_2^k$ with $|T|\leq 5^k$ such that $B_2^k\subset 2\conv (T)$.
\end{lem}

\begin{proof}
Let $T$ be the maximal $\frac{1}{2}$-separated set in $B_2^k$, the standard volumetric argument shows that
$|T|\leq 5^k$. We have $B_2^k\subset T+\frac{1}{2}B_2^k\subset \conv (T)+\frac{1}{2}B_2^k$, so
$B_2^k\subset 2\conv (T)$.
\end{proof}

The next lemma comes from \cite{Ko}.

\begin{lem}
\label{lem:smallconv}
For any $1\leq k\leq n$ there exists $T\subset B_2^n$ with $|T|\leq \frac{2n}{k}5^k$ such that
$B_2^n\subset 2\sqrt\frac{2n}{k} \conv(T)$.
\end{lem}

\begin{proof}
Let $l=\lceil n/k\rceil\leq 2n/k$ and $\er^n=F_1\oplus\cdots\oplus F_l$ be an orthogonal decomposition of $\er^n$ into
spaces of dimension at most $k$. By Lemma \ref{lem:net} we can find $T_i\subset B_2(F_i):= B_2^n\cap F_i$ such
that $B_2(F_i)\subset 2\conv(T_i)$ and $|T_i|\leq 5^k$. 
Let $T:=\bigcup_{i\leq l}T_i$. Then $T\subset B_2^n$ and $|T|\leq l5^k\leq  \frac{2n}{k}5^k$.

Fix now $x\in B_2^n$ and $x_i$ denotes its orthogonal projection on $F_i$. Observe that
\[
\sum_{i\leq l}\|x_i\|\leq \sqrt{l}\Big(\sum_{i\leq l}\|x_i\|^2\Big)^{1/2}\leq \sqrt{l}.
\]
Therefore
\[
x\subset \sqrt{l}\conv\Big\{0,\frac{x_1}{\|x_1\|},\ldots,\frac{x_l}{\|x_l\|}\Big\}
\subset \sqrt{l}\conv\Big(\bigcup_{i\leq l}B_2(F_i)\Big) \subset 2\sqrt{l}\conv(T).
\]
\end{proof}

\begin{lem}
\label{lem:expcutunif}
Let $Y$ be a vector uniformly distributed over $S^{n-1}$. Then 
\[
\Ex|\langle Y,t\rangle|I_{\{|\langle Y,t\rangle|\geq u\}}\leq
\min\Big\{\frac{|t|}{\sqrt{n}},\frac{2(|t|^2+nu^2)}{nu}e^{-nu^2/(2|t|^2)}\Big\}\quad t\in \er^n,\ u>0.
\]
\end{lem}

\begin{proof}
Observe that $\langle Y,t\rangle$ is distributed as $|t|Y_1$. Hence
\[
\Ex|\langle Y,t\rangle|I_{\{|\langle Y,t\rangle|\geq u\}}=|t|\Ex|Y_1|I_{\{|Y_1|\geq u/|t|\}}.
\]
We have $\Ex|Y_1|\leq (\Ex|Y_1|^2)^{1/2}=n^{-1/2}$. Moreover $\Pr(Y_1\geq v)\leq \exp(-nv^2/2)$ for $v\geq 0$
(cf. \cite{Tk}). Therefore
\begin{align*}
\Ex|Y_1|I_{\{|Y_1|\geq u\}}
&\leq u\Pr(|Y_1|\geq u)+\int_u^\infty\Pr(|Y_1|\geq v)\de v
\leq 2ue^{-nu^2/2}+ 2\int_u^\infty e^{-nv^2/2}\de v
\\
&\leq 2ue^{-nu^2/2}+2\int_u^\infty \frac{nv}{nu}e^{-nv^2/2}\de v
=\frac{2(1+nu^2)}{nu}e^{-nu^2/2}.
\end{align*}
\end{proof}

Now we are able to show that \eqref{eq:revest} holds for $T=B_2^n$ under $4+\delta$ moment condition.

\begin{prop}
\label{prop:b2n}
Let $X_1,\ldots,X_n$ be independent centered r.v's with variance 1 satisfying condition \eqref{eq:4+mom}. Then
there exists $S\subset \er^n$ such that $|S|\leq 10n^2$, $B_2^n\subset \conv(S)$ and
\[
M_X(S)\lesssim_{r,\lambda}\sqrt{n}\sim_{\lambda}\Ex|X|=b_X(B_2^n).
\] 
\end{prop}

\begin{proof}
By the Rosenthal inequality \cite{Ro} we have (recall that $r\in (4,8]$),
\begin{align*}
\big\||X|^2-n\big\|_{r/2}
&=\Bigg\|\sum_{i=1}^n(X_i^2-1)\Bigg\|_{r/2}\lesssim \Big(\sum_{i=1}^n\Var(X_i^2)\Big)^{1/2}+
\Big(\sum_{i=1}^n\Ex|X_i^2-1|^{r/2}\Big)^{2/r}
\\
&\lesssim_\lambda n^{1/2}+n^{2/r}\leq 2n^{1/2}.
\end{align*}
Therefore
\begin{equation}
\label{eq:expcutX}
\Ex|X|I_{\{|X|\geq \sqrt{2n}\}}\leq \Ex\sqrt{2(|X|^2-n)}I_{\{|X|\geq \sqrt{2n}\}}
\leq \sqrt{2}n^{1/2-r/2}\Ex(|X|^2-n)^{r/2}\leq C(\lambda)n^{1/2-r/4}.
\end{equation}

By Lemma \ref{lem:smallconv} (applied with $k=c(r)\log n$) there exists $t_1,\ldots,t_N$  such that
$B_2^n\subset \conv\{t_1,\ldots,t_N\}$, $N\leq 10n^{1/2+r/8}$ and $|t_i|\leq C(r)\sqrt{n/\log n}$, $1\leq i\leq N$. Let $U$ be the random rotation
(uniformly distributed on $O(n)$) then $Ut_i$ is distributed as $|t_i|Y$, where $Y$ has uniform distribution on $S^{n-1}$.
Thus by Lemma \ref{lem:expcutunif},
\begin{align*}
\Ex_U\Ex_X|\langle X,Ut_i\rangle|I_{\{|\langle X,Ut_i\rangle|\geq u\}}
&=\Ex_X\Ex_Y|\langle Y,|t_i|X\rangle|I_{\{|\langle Y,|t_i|X\rangle|\geq u\}}
\\
&\leq \Ex\min\Big\{\frac{|t_i||X|}{\sqrt{n}},\frac{2(|t_i|^2|X|^2+nu^2)}{nu}e^{-nu^2/(2|t_i|^2|X|^2)}\Big\}
\\
&\leq \frac{|t_i|}{\sqrt{n}}\Ex|X|I_{\{|X|\geq \sqrt{2n}\}}+\frac{4|t_i|^2+2u^2}{u}e^{-u^2/(4|t_i|^2)}.
\end{align*}

Recall that $|t_i|\lesssim_r \sqrt{n/\log n}$ so for sufficiently large $C(r)$ we get by \eqref{eq:expcutX},
\[
\Ex_U\Ex_X|\langle X,Ut_i\rangle|I_{\{|\langle X,Ut_i\rangle|\geq C(r)\sqrt{n}\}}
\leq C(\lambda)n^{-r/4}|t_i|+n^{-2}\leq C(r,\lambda)n^{1/2-r/4}.
\]

As a consequence there exists $U\in O(n)$ such that 
\begin{equation}
\label{eq:estMXball}
\sum_{i=1}^N \Ex_X|\langle X,Ut_i\rangle|I_{\{|\langle X,Ut_i\rangle|\geq C(r)\sqrt{n}\}}
\leq NC(r,\lambda)n^{1/2-r/4}\leq 10C(r,\lambda)n^{1-r/8}.
\end{equation}
Thus if we put $S:=\{Ut_1,\ldots,Ut_N\}$ we will have $\conv(S)=U\conv\{t_1,\ldots,t_N\}\supset B_2^n$ and
$M_X(S)\leq C'(r,\lambda)\sqrt{n}$.

\end{proof}

\subsection{Ellipsoids}

We now extend the bounds from the previous subsection to the case of ellipsoids, i.e. sets of the form
\begin{equation}
\label{eq:defell}
\mathcal{E}:=\Big\{t\in \er^n\colon\ \sum_{i=1}^n\frac{\langle t,u_i\rangle^2 }{a_i^2}\leq 1\Big\},
\end{equation}
where $u_1,\ldots,u_n$ is an orthonormal system in $\er^n$ and $a_1,\ldots,a_n>0$.

Observe that 
\[
\sup_{t\in \mathcal{E}}\langle t,x\rangle=\sqrt{\sum_{i=1}^n a_i^2\langle x,u_i\rangle ^2}.
\]

To treat this case we will need the following Lemma.

\begin{lem}
\label{lem:ellips1}
Let $X=(X_1,\ldots,X_n)$, where $X_i$ are independent mean zero and variance one r.v's satisfying $4+\delta$ condition \eqref{eq:4+mom}. \\
i) For any $a_1,\ldots,a_n\geq 0$ and any o.n. vectors $u_1,\ldots,u_n$,
\[
\Ex\Big(\sum_{k=1}^n a_k^2 \langle X,u_k\rangle^2\Big)^{1/2}\sim_{\lambda}
\Big(\Ex\sum_{k=1}^n a_k^2  \langle X,u_k\rangle^2\Big)^{1/2}=\Big(\sum_{k=1}^n a_k^2\Big)^{1/2}.
\]
ii) For any $n\times n$ matrix $B$,
\[
\Big(\Ex(|BX|^2-\|B\|_{\mathrm{HS}}^2)^{r/2}\Big)^{2/r}\leq C(\lambda) \|B^TB\|^{1/2}_{\mathrm{HS}}.
\]
In particular for any linear supspace $E\subset \er^n$ od dimension $k\in \{1,\ldots,n\}$,
\[
\Big(\Ex(|P_EX|^2-k)^{r/2}\Big)^{2/r}\leq C(\lambda) k^{1/2}.
\]
\end{lem}

\begin{proof} Part i) follows from Lemma \ref{lem:estmom4+}.

To show part ii) let $B=(b_{ij})_{i,j=1}^n$, $e_1,e_2,\ldots,e_n$ be the canonical basis of $\er^n$ 
and let
\[
\sigma_{i,j}:=\sum_{l=1}^nb_{l,i}b_{l,j}=\langle Be_i,Be_j\rangle,\quad 1\leq i,j\leq n.
\]
Then
\begin{align*}
\big\||BX|^2-\|B\|_{\mathrm{HS}}^2\big\|_{r/2}
&=\Big\|\sum_{i=1}^n(X_i^2-1)\sigma_{i,i}+\sum_{1\leq i\neq j\leq n}X_{i}X_j\sigma_{i,j}\Big\|_{r/2}
\\
&\leq \Big\|\sum_{i=1}^n(X_i^2-1)\sigma_{i,i}\Big\|_{r/2}+\Big\|\sum_{1\leq i\neq j\leq n}X_{i}X_j\sigma_{i,j}\Big\|_{r/2}.
\end{align*}

Applying Rosenthal's inequality we get
\begin{align*}
\Big\|\sum_{i=1}^n(X_i^2-1)\sigma_{i,i}\Big\|_{r/2}
&\lesssim \Big(\sum_{i=1}^n\Var(X_i^2)\sigma_{i,i}^2\Big)^{1/2}+\Big(\sum_{i=1}^n\Ex(X_i^2-1)^{r/2}\sigma_{i,i}^{r/2}\Big)^{2/r}
\\
&\lesssim_\lambda \Big(\sum_{i=1}^n\sigma_{i,i}^2\Big)^{1/2}+\Big(\sum_{i=1}^n\sigma_{i,i}^{r/2}\Big)^{2/r}\leq 
2\Big(\sum_{i=1}^n\sigma_{i,i}^2\Big)^{1/2}.
\end{align*}

Hypercontractive method (as in the proof of Lemma \ref{lem:estmom4+}) yields
\[
\Big\|\sum_{i\neq j}X_{i}X_j\sigma_{i,j}\Big\|_{r/2}
\lesssim_\lambda \Big\|\sum_{i\neq j}X_{i}X_j\sigma_{i,j}\Big\|_{2}=
\Big(\sum_{i\neq j}\sigma_{i,j}^2\Big)^{1/2}.
\]
Finally 
\[
\Big(\sum_{i=1}^n\sigma_{i,i}^2\Big)^{1/2}+\Big(\sum_{i\neq j}\sigma_{i,j}^2\Big)^{1/2}\leq 
2\Big(\sum_{i,j}\sigma_{i,j}^2\Big)^{1/2}=2\|B^TB\|_{\mathrm{HS}}.
\]
\end{proof}

Now we state and prove the main result of this section.

\begin{thm}
\label{thm:ellipsoid}
Let $X_1,\ldots,X_n$ be independent centered r.v's satisfying the condition \eqref{eq:4+mom} and let $T$ be an ellipsoid in
$\er^n$. 
Then there exists $S\subset \er^n$ such that $|S|\leq 10n^2$, $T\subset \conv(S)$ and
\[
M_X(S)\lesssim_{r,\lambda}b_X(T).
\] 
\end{thm}

\begin{proof} Since it is only a matter of scaling we may and will assume that $\Ex X_i^2=1$ for all $i$. 
Let $T=\mathcal{E}$ be an ellipsoid of the form \eqref{eq:defell}. Then the first part of Lemma \ref{lem:ellips1}
yields
\[
\Ex \sup_{t\in \mathcal{E}}X_t=\Ex\Big(\sum_{k=1}^n a_k^2 \langle X,u_k\rangle^2\Big)^{1/2}\sim_{\lambda}
\sqrt{\sum_{k=1}^n a_k^2}.
\]
By homogenity we may assume that $\sum_{k=1}^n a_k^2=1$.

Define 
\[
I_k:=\{i\colon\ 2^{-k-1}< a_i\leq 2^{-k}\},\ \ n_k:=|I_k|,\ \ J:=\{k\in\zet\colon\ I_k\neq \emptyset\}, \ \
E_k:=\mathrm{span}\{u_i\colon i\in I_k\}.
\]
Then 
\begin{equation}
\label{eq:norm2}
1\leq \sum_{k\in J}n_k2^{-2k}<4.
\end{equation}
In particular $J$ is a subset of nonnegative integers.

We claim that for any positive sequence $(c_k)_{k\in J}$ such that $\sum_{k}c_{k}^{-2}\leq 1$,
\[
\mathcal{E}\subset \conv\Big(\bigcup_{k\in J}c_k2^{-k}B_2^{I_k}\Big), \quad \mbox{where}\quad B_2^{I_k}:=B_2^n\cap E_k. 
\]
Indeed, let $P_kx:=\sum_{i\in I_k}\langle x,u_i\rangle u_i$ be the projection of $x$ onto $E_k$, then 
\[
x=\sum_{k\in J}c_{k}^{-1}2^k|P_kx|c_k2^{-k}\frac{P_kx}{|P_kx|}
\]
and for $x\in \mathcal{E}$,
\[
\sum_{k\in J}c_{k}^{-1}2^k|P_kx|\leq \sqrt{\sum_{k\in J}c_{k}^{-2}}\sqrt{\sum_{k\in J}2^{2k}|P_kx|^2}\leq
\sqrt{\sum_{k\in J}\sum_{i\in I_k}\frac{\langle x,u_i\rangle^2}{a_i^2}}\leq 1.
\]

Let us for a moment fix $k\in J$. By Lemma \ref{lem:smallconv} (applied with $k=c(r)\log n_k$) there exists 
$t_1,\ldots,t_{N_k}\in E_k$  such that
$B_2^{I_k}\subset \conv\{t_1,\ldots,t_{N_k}\}$, $N_k\leq 10n_k^{1/2+r/8}$ and $|t_i|\leq C(r)\sqrt{n_k/\log (n_k)}$. 
Let $U$ be the random rotation of $E_k$ (uniformly distributed on $O(E_k)$) 
then $Ut_i$ is distributed as $|t_i|Y$, where $Y$ has uniform distribution on 
$S^{I_k}:=S^{n-1}\cap E_k$.
Thus by Lemma \ref{lem:expcutunif},
\begin{align*}
\Ex_U\Ex_X|\langle X,Ut_i\rangle|&I_{\{|\langle X,Ut_i\rangle|\geq u\}}
\\
&=\Ex_X\Ex_Y|\langle Y,|t_i|P_{E_k}X\rangle|I_{\{|\langle Y,|t_i|P_{E_k}X\rangle|\geq u\}}
\\
&\leq \Ex\min\Big\{\frac{|t_i||P_{E_k}X|}{\sqrt{n_k}},\frac{2(|t_i|^2|P_{E_k}X|^2+n_ku^2)}{n_ku}e^{-n_ku^2/(2|t_i|^2|P_{E_k}X|^2)}\Big\}
\\
&\leq \frac{|t_i|}{\sqrt{n_k}}\Ex|P_{E_k}X|I_{\{|P_{E_k}X|\geq \sqrt{2n_k}\}}+\frac{4|t_i|^2+2u^2}{u}e^{-u^2/(4|t_i|^2)}.
\end{align*}
We have
\begin{align*}
\Ex|P_{E_k}X|I_{\{|P_{E_k}X|\geq \sqrt{2n_k}\}}
&\leq \sqrt{2}\Ex(|P_{E_k}X|^2-n_k)^{1/2}I_{\{|P_{E_k}X|\geq \sqrt{2n_k}\}}
\\
&\leq \sqrt{2}n_k^{1/2-r/2}\Ex(|P_{E_k}X|^2-n_k)^{r/2}
\leq C(\lambda)n_k^{1/2-r/4},
\end{align*}
where the last bound follows by Lemma \ref{lem:ellips1}.
Recall that $|t_i|\lesssim_r \sqrt{n_k/\log n_k}$, thus for sufficiently large $C(r)$ we get 
\[
\Ex_U\Ex_X|\langle X,Ut_i\rangle|I_{\{|\langle X,Ut_i\rangle|\geq C(r)\sqrt{n}\}}
\leq C(\lambda)n_k^{-r/4}|t_i|+n_k^{-2}\leq C(r,\lambda)n_k^{1/2-r/4}.
\]

As a consequence there exists $U\in O(E_k)$ such that 
\[
\sum_{i=1}^{N_k} \Ex_X|\langle X,Ut_i\rangle|I_{\{|\langle X,Ut_i\rangle|\geq C(r)\sqrt{n_k}\}}
\leq N_kC(r,\lambda)n_k^{1/2-r/4}\leq 10C(r,\lambda)n_k^{1-r/8}.
\]
Define $S_k=\{t_{k,1},\ldots,t_{k,N_k}\}:=\{Ut_1,\ldots,Ut_{N_k}\}$. Then 
$\conv(S_k)=U\conv\{t_1,\ldots,t_{N_k}\}\supset B_2^{I_k}$, $N_k\leq 10n_k^{1/2+r/8}\leq 10n_k^2$  and
\[
\sum_{i=1}^{N_k} \Ex_X|\langle X,t_{k,i}\rangle|I_{\{|\langle X,t_{k,i}\rangle|\geq C(r)\sqrt{n_k}\}}
\leq C(r,\lambda)n_k^{1-r/8}.
\]

Set $c_k:=2^{k+2}(2^k+n_k)^{-1/2}$. By \eqref{eq:norm2} we get $\sum_{k\in J}c_{k}^{-2}\leq 1$, so
\[
\mathcal{E}
\subset \conv\Big(\bigcup_{k\in J}c_k2^{-k}B_2^{I_k}\Big)\subset \conv(\{c_k2^{-k}t_{k,i}\colon\ k\in J,i\leq N_k\}):=\conv(S).
\]
We have 
\[
|S|= \sum_{k\in J}N_k\leq \sum_{k\in J}10n_k^2\leq 10\Big(\sum_{k\in J}n_k\Big)^2=10n^2.
\]
Moreover,
\begin{align*}
\sum_{s\in S}\Ex|\langle s,X\rangle|I_{\{|\langle s,X\rangle|\geq 4C(r)\}}
&
=\sum_{k\in J}2^{-k}c_k\sum_{i=1}^{N_k}\Ex|\langle t_{k,i},X\rangle|I_{\{2^{-k}c_k|\langle t_{k,i},X\rangle|\geq 4C(r)\}}
\\
&\leq
\sum_{k\in J}4(2^k+n_k)^{-1/2}\sum_{i=1}^{N_k}\Ex|\langle t_{k,i},X\rangle|I_{\{|\langle t_{k,i},X\rangle|\geq C(r)\sqrt{n_k}\}}
\\
&\leq \sum_{k\in J}4(2^k+n_k)^{-1/2}C(r,\lambda)n_k^{1-r/8}
\\
&\leq 4C(r,\lambda)\sum_{k\in J}(2^k+n_k)^{1/2-r/8}
\leq 4C(r,\lambda)\sum_{k\geq 0}2^{k(1/2-r/8)}
\\
&\leq C'(r,\lambda),
\end{align*}
which shows that $M_X(S)\sim \widetilde{M}_X(S)\lesssim_{\lambda,r}1\sim b_X(\mathcal{E})$.

\end{proof}

\section{Case III.  $\ell_q^n$-balls, $2<q\leq \infty$}
\label{sec:Bq}

It turns out that results of the previous sections may be easily applied to get estimates in the case when $T=B_q^n$
is the unit ball in $\ell_q^n$ and $q\in (2,\infty]$. In the whole section by $q'$ we will denote the H\"older dual of $q$, i.e.
$q'=\frac{q}{q-1}$, $2\leq q<\infty$ and $q'=1$ for $q=\infty$.

\begin{prop}
\label{prop:bqn}
Let $X_1,\ldots,X_n$ be independent centered r.v's with variance 1 satisfying condition \eqref{eq:4+mom}. Then
there exists $S\subset \er^n$ such that $|S|\leq 10n^2$, $B_q^n\subset \conv(S)$ and
\[
M_X(S)\lesssim_{r,\lambda}n^{1/q'}\sim_\lambda b_X(B_q^n).
\] 
\end{prop}

\begin{proof}
Since $q'\in (1,2]$, condition \eqref{eq:4+mom} yields $\|X_i\|_{q'}\sim_\lambda \|X_i\|_{2q'}\sim_\lambda\|X_i\|_2=1$ and hence
$(\Ex\|X\|_{q'}^{2q'})^{1/(2q')}\sim_\lambda (\Ex\|X\|_{q'}^{q'})^{1/q'}$. Therefore
\[
b_X(B_q^n)=\Ex\sup_{t\in B_q^n}\langle t,X\rangle =\Ex\|X\|_{q'}\sim_\lambda \big(\Ex\|X\|_{q'}^{q'}\big)^{1/q'}\sim_\lambda n^{1/q'}.
\]

H\"older's inequality implies $B_q^n\subset n^{1/2-1/q}B_2^{n}=n^{1/q'-1/2}B_2^n$ and the assertion easily follows from Proposition \ref{prop:b2n}.
\end{proof}

Now let us consider the case of linear transformation of $\ell_q^n$-ball, i.e. $T=AB_q^n$.  Next simple lemma shows how to estimate $b_X(T)$.

\begin{lem}
\label{lem:bXABqn}
Let $X=(X_1,\ldots,X_n)$, where $X_i$ are independent mean zero and variance one r.v's satisfying $4+\delta$ condition \eqref{eq:4+mom}. Then for any $n\times n$ matrix $A$ and $2\leq q\leq \infty$ we have
\[
b_X(AB_q^n)=b_{A^TX}(B_q^n)\sim_{\lambda}\Big(\sum_{i=1}^n |Ae_i|^{q'}\Big)^{1/q'}.
\]
\end{lem}

\begin{proof}
Observe that
\[
\sup_{t\in AB_q^n}\langle X,t\rangle=\sup_{t\in B_q^n}\langle A^TX,t\rangle
=\Big(\sum_{i=1}^n|\langle A^TX,e_i\rangle |^{q'}\Big)^{1/q'}
=\Big(\sum_{i=1}^n|\langle X,Ae_i\rangle |^{q'}\Big)^{1/q'}.
\]
Condition \eqref{eq:4+mom} (see Lemma \ref{lem:estmom4+}) implies that 
\[
\|\langle X,Ae_i\rangle \|_{2q'}\sim_\lambda \|\langle X,Ae_i\rangle \|_{q'}\sim_\lambda \|\langle X,Ae_i\rangle \|_{2}
=|Ae_i|.
\]
Hence $\|\sup_{t\in AB_q^n}\langle X,t\rangle\|_{2q'}\sim_\lambda \|\sup_{t\in AB_q^n}\langle X,t\rangle\|_{q'}$
and 
\[
b_X(AB_q^n)=\Big\|\sup_{t\in AB_q^n}\langle X,t\rangle\Big\|_{1}
\sim_\lambda \Big\|\sup_{t\in AB_q^n}\langle X,t\rangle\Big\|_{q'}
\sim_\lambda \Big(\sum_{i=1}^n |Ae_i|^{q'}\Big)^{1/q'}.
\]

\end{proof}

As in the proof of Proposition \ref{prop:bqn} we may include linear image of $B_q^n$ into ellipsoid with the comparable $b_X$-bound and deduce from Theorem \ref{thm:ellipsoid} the following more general result. 

\begin{thm}
Let $X_1,\ldots,X_n$ be independent centered r.v's satisfying condition \eqref{eq:4+mom} and let $T=AB_q^n$ for some 
$2\leq q\leq \infty$ and an $n\times n$ matrix $A$. 
Then there exists $S\subset \er^n$ such that $|S|\leq 10n^2$, $T\subset \conv(S)$ and
\[
M_X(S)\lesssim_{r,\lambda}b_X(T).
\] 
\end{thm}

\begin{proof} Since it is only a matter of scaling we may and will assume that $\Ex X_i^2=1$ for all $i$. 
By Lemma \ref{lem:bXABqn} it is enough to show that
\[
M_X(S)\lesssim_{r,\lambda}\Big(\sum_{i=1}^n |Ae_i|^{q'}\Big)^{1/q'}.
\] 
By homogenity we may assume that $\sum_{i=1}^n |Ae_i|^{q'}=1$. Case $q=2$ was treated in Theorem  \ref{thm:ellipsoid}, so
we may assume that $q>2$, i.e. $q'<2$. Moreover, we may assume that $Ae_i\neq 0$ for all $i$.

Let  $\lambda_i:=|Ae_i|^{1-q'/2}$. Observe that if $t\in B_q^n$ then by H\"older's inequality
\begin{align*}
\sum_{i=1}^n|\lambda_it_i|^2
&\leq \Big(\sum_{i=1}^n |t_i|^q\Big)^{2/q}\Big(\sum_{i=1}^n |\lambda_i|^{2q/(q-2)}\Big)^{(q-2)/q}
\\
&=\Big(\sum_{i=1}^n |t_i|^q\Big)^{2/q}\Big(\sum_{i=1}^n |Ae_i|^{q'}\Big)^{(q-2)/q}\leq 1.
\end{align*}
This shows that $D^{-1}B_q^n\subset B_2^n$, where $D:=\mathrm{diag}(d_1,\ldots,d_n)$ and $d_i:=|Ae_i|^{q'/2-1}$. 
Hence  $AB_q^n\subset ADB_2^n$ and
\[
b_X(ADB_2^n)\sim_{\lambda} \Big(\sum_{i=1}^n|ADe_i|^2\Big)^{1/2}=\Big(\sum_{i=1}^n|Ae_i|^{q'}\Big)^{1/2}=1.
\]
We get the assertion applying Theorem \ref{thm:ellipsoid} for the ellipsoid $ADB_2^n$.

\end{proof}

\section{Concluding remarks and open questions}

We have shown that the main question has the affirmative answer in the case $T$ is an ellipsoid (or more general linear image of $\ell_q^n$-ball, $2\leq q\leq n$) if $X_i$ are independent mean zero r.v's satisfying the
$4+\delta$ moment condition \eqref{eq:4+mom}. The following questions are up to our best knowledge open.

\begin{itemize}
\item Does \eqref{eq:revest} holds for $T=B_q^n$, $1<q<2$ and $X_i$ satisfying $4+\delta$ moment condition?

\item John's theorem states that for any convex symmetric set $T$ in $\er^n$ there exists an ellipsoid $\mathcal{E}$
such that $\mathcal{E}\subset T\subset \sqrt{n}\mathcal{E}$. Hence Theorem \ref{thm:ellipsoid} implies that under
$4+\delta$ condition \eqref{eq:4+mom} one may
find finite set $S$ such that $T\subset \conv(S\cup -S)$ and $M_X(S)\leq C(r,\lambda)\sqrt{n}b_X(T)$.
We do not whether one may improve upon  $\sqrt{n}$ factor for general sets $T$.

\item Are there heavy-tailed random variables $X_i$ such that \eqref{eq:revest} holds for arbitrary set $T$
(for heavy-tailed r.v's approach via chaining functionals described in Subsection \ref{sub:gammaX} fails to work)? 

\item Let $X_i$ be heavy-tailed symmetric Weibull r.v's (i.e. symmetric variables with tails $\exp(-t^r)$,
$0<r<1$). Bogucki \cite{Bo} was able to obtain two-sided bounds for $b_X(T)$ with the use of random permutations (which may be eliminated if $T$ is permutationally invariant). We do not know if the convex hull method works in this case.
 
\end{itemize}

\end{document}